\documentclass[12pt, reqno]{amsart}
\usepackage{amsmath, amsthm, amscd, amsfonts, amssymb, graphicx, color}
\usepackage[bookmarksnumbered, colorlinks, plainpages]{hyperref}
\usepackage{comment}

\usepackage{hyperref}
\usepackage{longtable}
\usepackage{xcolor,soul}
\usepackage{listings}

\definecolor{dkgreen}{rgb}{0,0.6,0}
\definecolor{gray}{rgb}{0.5,0.5,0.5}
\definecolor{mauve}{rgb}{0.58,0,0.82}

\newtheorem{theorem}{Theorem}[section]
\newtheorem{lemma}[theorem]{Lemma}

\theoremstyle{definition}
\newtheorem{definition}[theorem]{Definition}

\theoremstyle{remark}
\newtheorem{remark}[theorem]{Remark}
\numberwithin{equation}{section}

\begin{document}
\setcounter{page}{1}

\noindent {\small Research Preprint}\hfill     {\small ISSN: XXXX-YYYY}\\
{\small Vol XX, Issue Y(202Z) 1-19}\hfill  {\small https://doi.org/}

\centerline{}

\centerline{}

\title[EDGE-GRACEFUL USUAL FAN GRAPHS]{EDGE-GRACEFUL USUAL FAN GRAPHS}

\author[A. Angel, J. Antalan, J. Gamurot, R. Tagle]{Aaron D.C. Angel$^1$, John Rafael M. Antalan$^2$$^{*}$, John Loureynz F.Gamurot$^1$, \MakeLowercase {and} Richard P. Tagle$^2$}

%In case of 3 or more authors use below format
%\author[F. Author, S. Author, T. Author]{First Author$^1$, Second Author$^2$$^{*}$ \MakeLowercase {and} Third Author$^3$}

\address{$^{1}$ Alumni, Department of Mathematics and Physics, College of Science, Central Luzon State University (3120), Science City of Mu\~{n}oz, Nueva Ecija, Philippines.}
\email{\textcolor[rgb]{0.00,0.00,0.84}{aaronangel1123@gmail.com,angel.aaron@clsu2.edu.ph}}
\email{\textcolor[rgb]{0.00,0.00,0.84}
{gamurotloureynz@gmail.com, gamurot.john@clsu2.edu.ph}}

\address{$^{2}$ Faculty, Department of Mathematics and Physics, College of Science, Central Luzon State University (3120), Science City of Mu\~{n}oz, Nueva Ecija, Philippines.} 
\email{\textcolor[rgb]{0.00,0.00,0.84}{jrantalan@clsu.edu.ph}}
\email{\textcolor[rgb]{0.00,0.00,0.84}{richard\_tagle@clsu.edu.ph}}

%\dedicatory{This paper is dedicated to Professor ABCD}



\date{Received: xxxxxx; Revised: yyyyyy; Accepted: zzzzzz.
\newline \indent $^{*}$ Corresponding author}

\begin{abstract}
A graph $G$ with $p$ vertices and $q$ edges is said to be edge-graceful if its edges can be labeled from $1$ through $q$, in such a way that the labels induced on the vertices by adding over the labels of incident edges modulo $p$ are distinct. A known result under this topic is Lo's Theorem, which states that if a graph $G$ with $p$ vertices and $q$ edges is edge-graceful, then $p\Big|\Big(q^{2}+q-\dfrac{p(p-1)}{2}\Big)$.

This paper presents novel results on the edge-gracefulness of the usual fan graphs. Using Lo's Theorem, the concepts of divisibility and Diophantine equations, and a computer program created, we determine all edge-graceful usual fan graphs $F_{1,n}$ with their corresponding edge-graceful labels. 
\end{abstract} \maketitle

\section{Introduction}

The study of graphs is a vast field in mathematics which offers a plethora of unique and interesting topics. One of these particular topics in graph theory is graph labeling.

The concept of {\bf graph labeling} was first introduced in the mid-1960s. Graph labeling is an assignment of integers to the vertices or edges, or both, of a graph under certain conditions \cite{gallian}. In the last 60 years, over 200 types or variations of graph labeling have been studied and about 2500 articles have been published \cite{rama}. One type of graph labeling is graceful labeling. A {\bf graceful labeling} or graceful numbering is a special graph labeling of a graph on $m$ edges in which the vertices are labeled, using a subset of distinct nonnegative integers from $0$ to $m$ and each edge of the graph is uniquely labeled by the absolute difference between the labels of the vertices incident to it. If the resulting graph edge numbers run from $1$ to $m$ inclusive, it is a graceful labeling, and the graph is said to be a graceful graph \cite{ven}.

On the other hand, the concept of {\bf edge-gracefulness} and edge-graceful graphs was defined and introduced by Lo in 1985 \cite{Lo}. Edge-graceful labeling is considered as a reversal of graceful labeling because it labels the edges first; then the labels of the vertices would depend on the labels of the edges incident to them. That is, a graph $G$ with $p$ vertices and $q$ edges is said to be edge-graceful if the edges can be labeled from $1$ through $q$, in such a way that the labels induced on the vertices by adding over the labels of incident edges modulo $p$ are distinct \cite{ven}.

Since the introduction of the concept of edge-graceful labeling in 1985, several works on the edge-graceful labeling of specific families of graphs have been conducted. For instance, Kuan et al. \cite{kuan} studied edge-graceful unicyclic graphs. Lee et al. \cite{lee} on the other hand, studied the edge-graceful labeling of complete graphs. The 
work of Small in \cite{small} shows that regular (even) spider graphs are edge-graceful. In 1992, Cabannis et al. \cite{caban} studied the edge-gracefulness of regular graphs and trees. Fast-forward to 2017, Venkatesan and Sekar \cite{ven} attempted to prove that among all wheel graphs $W_n$, only $W_3$ is edge-graceful while a complete proof of their result was provided recently in \cite{angel}.    

In this paper, we used essential graph theory concepts to study the basic properties of usual fan graphs.  In particular, their vertex and edge sets were used to apply Lo's Theorem, an initial condition that must be satisfied to conclude that a particular graph is edge-graceful.  Then, using the definition of edge-graceful labeling, and some preliminary concepts in number theory, which include the concept of the Diophantine equation, this paper determined which among the said graphs satisfied Lo's Theorem. If a particular graph does not satisfy Lo's Theorem, then it automatically does not qualify to be an edge-graceful graph. Finally, by using a computer program that we created, this paper constructed actual examples of edge-graceful labeling of the usual fan graphs that satisfy these initial conditions, thereby showing that the graphs are indeed edge-graceful.

\section{Preliminaries}

\label{sec2}

For completeness, we provide all the preliminary concepts and results necessary for the clear and organized presentation of the results of this study. Moreover, we emphasize that the graph-theoretical definitions used in this section, aside from the concept of edge-graceful labeling, are from Wolfram MathWorld \cite{Wolfram}.
We also refer the readers who are unfamiliar with graph theory concepts that were briefly stated in this paper, to some standard graph theory references such as Gallian's ``A Dynamic Survey of Graph Labeling" \cite{gallian} and Gross's ``Graph theory and its applications"\cite{gross}.

\subsection{Essential Graph Theory Concepts}
\label{subsec2}

To begin this subsection, we briefly discuss some key concepts in graph theory that are essential to this paper.

\begin{definition}
	A $(p,q)$-graph $G(V,E)$ is called \textbf{edge-graceful} if there exist a bijection $f:E\rightarrow \{1,2,....,q\}$, such that the induced mapping $f^{+} :V \rightarrow \mathbb{Z}_{p}$, defined by
	\begin{equation*}
		f^{+}(u)\equiv \Big(\sum_{v\in N(u)}f(uv)\Big)(mod\hspace{1mm}p)
	\end{equation*}
	is also a bijection, where $\mathbb{Z}_{p}$ is the set of integers $\{0,1,2,..,(p-1)\}$. Moreover, every $f(uv)$ are the labels of the edges with endpoint $u$ in $G$.
\end{definition}

\begin{figure}[h]
	\centering
	\includegraphics[width=8cm, height=4cm]{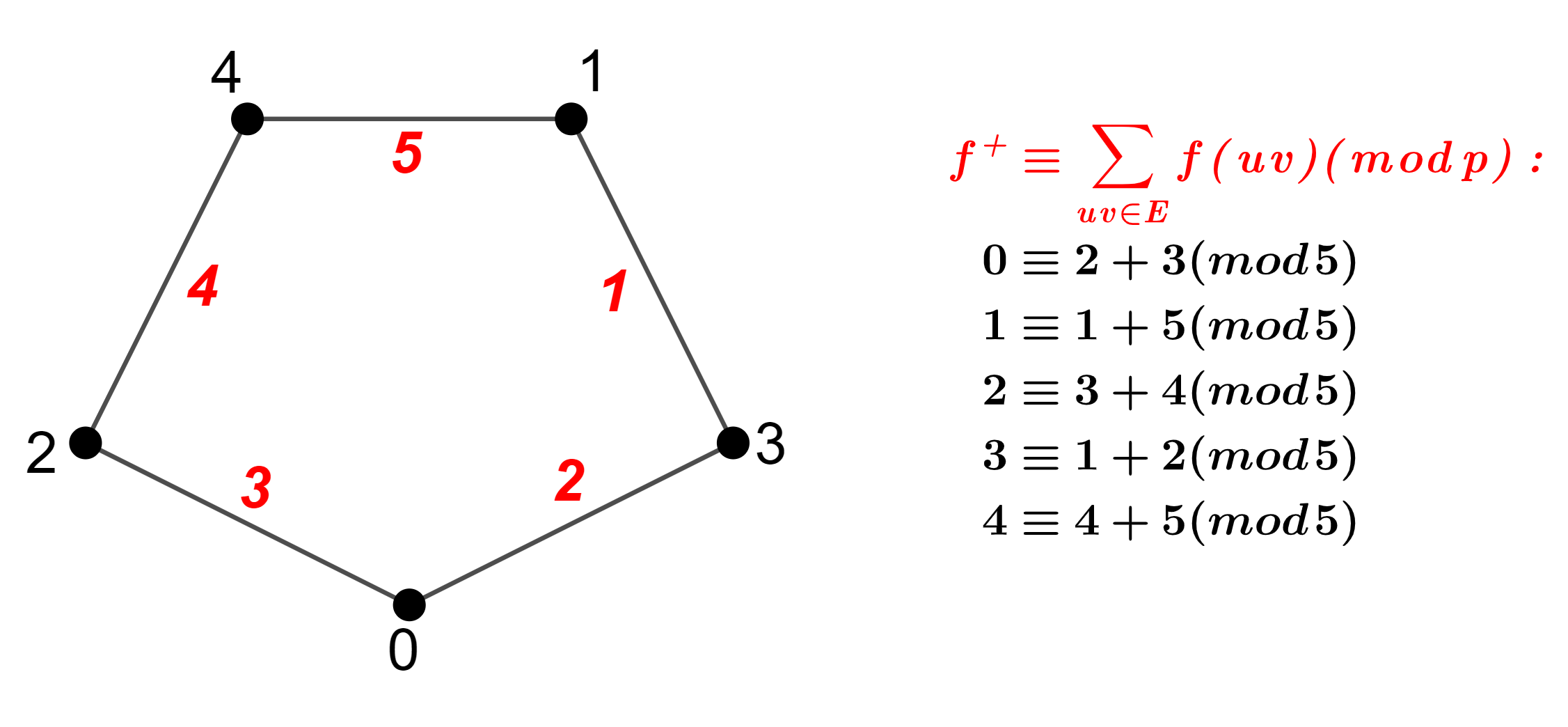}
	\caption{An edge-graceful labeling of the cycle graph $C_{5}$}\label{fig1.1}
\end{figure}

The {\bf{Figure \ref{fig1.1}}} in the next page, shows an example of an edge-graceful labeling of a graph. First, the graph involved, has five edges and five vertices. By following  {\bf{Definition 2.1}}, the edges of the graph were labeled using distinct integers from 1 to 5, and the induced labels of the vertices are distinct integers from 0 to 4.

\begin{definition}
	A \textbf{fan graph} denoted by $F_{m,n}$ is defined as the graph join $K_m+P_n$, where $K_m$ is the empty graph on $m$ vertices and $P_n$ is the path graph on $n$ vertices. The case when $m = 1$ corresponds to the usual fan graphs.
\end{definition}

Some examples of usual fan graphs are shown in {\bf{Figure \ref{fig2}}} below.

\begin{figure}[h]
	\centering
	\includegraphics[width=13cm, height=3.5cm]{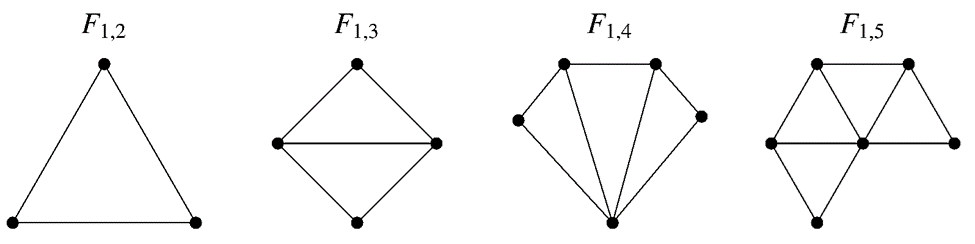}
	\caption{Examples of usual fan graphs}\label{fig2}
\end{figure}

\subsection{Essential Number Theory Concepts}
\label{subsec2}

For this subsection, we enumerate the important concepts in number theory that were utilized in this paper. For a more detailed discussion of basic concepts in number theory, readers may refer to Burton's ``Elementary Number Theory" \cite{burton}.

\begin{definition}
	Let $a$ and $b$ be two positive integers, where $a\leq b$. We say that $a$ {\bf divides} $b$, written in symbols by $a\mid b$, if there is a positive integer $c$ such that $b=ac$. Otherwise, we say that $a$ {\bf does not divide} $b$, denoted by $a\nmid b$. Now, if $a\mid b$, it also implies the following: (i) $b$ is a multiple of $a$, (ii) $a$ is a divisor of $b$ , and (iii) $a$ is a factor of $b$.
\end{definition}

\begin{definition}
	Let $p$ be a positive integer. Also, let $a,b\in \mathbb{Z}$. We can say that $a$ is congruent to $b$ (modulo $p$) or $a\equiv b\pmod{p}$ if $p\mid (a-b)$. That is, if $a=b+kp$ for some $k\in \mathbb{Z}$.
\end{definition}

\begin{theorem} \label{thm6}
	Let $a,b\in \mathbb{Z}$, and $p$ be a positive integer. If $a\equiv b(mod\hspace{1mm}p)$, then $b\equiv a(mod\hspace{1mm} p)$.
\end{theorem}

\begin{theorem}\label{thm7}
	Suppose $a\equiv b(mod\hspace{1mm}p)$ and $c\equiv d(mod\hspace{1mm}p)$, then $a+c\equiv b+d(mod\hspace{1mm}p)$
\end{theorem}

\begin{definition}
	{\bf General Quadratic Diophantine Equation}: Diophantine equations are equations wherein only integer solutions are allowed. A general quadratic Diophantine equation is an equation of the form $ax^{2} + bxy + cy^{2} + dx + ey + f = 0$, where both $x$ and $y$ are integers.
\end{definition}

\subsection{Essential Results Used in this Paper}
\label{subsec3}
This subsection discusses and proves Lo’s theorem, which is the main theorem that we will use to obtain the results of this paper. Furthermore, some essential transformations of quadratic Diophantine equations are also discussed in this section.

\begin{theorem}[Lo's Theorem\cite{Lo}]\label{thm9}
	If a graph $G$ of $p$ vertices and $q$ edges is edge-graceful, then
	\begin{equation*}
		p\Big|\Big(q^{2}+q-\dfrac{p(p-1)}{2}\Big).
	\end{equation*} 
\end{theorem}

\noindent \textit{Proof.} Again, we note that a $(p,q)$-graph $G(V,E)$ is called edge-graceful if there exist a bijection $f:E\rightarrow \{1,2,....,q\}$, such that the induced mapping $f^{+} :V \rightarrow \mathbb{Z}_{p}$, defined by
\begin{equation*}
	f^{+}(u)\equiv\Big( \sum_{v\in N(u)}f(uv)\Big)(mod\hspace{1mm}p)
\end{equation*}
is also a bijection, where every $f(uv)$ are the labels of the edges with endpoint $u$ in $G$.\\
Now, we denote by $u_{1}, u_{2},..., u_{p}$, the vertices in our graph with labels
\begin{equation*}
	f^{+}(u_{i})\equiv\Big( \sum_{v\in N(u_{i})}f(u_{i}v)\Big)(mod\hspace{1mm}p)\hspace{3mm},\hspace{1mm}i=1,2,...,p.
\end{equation*}
Adding all the labels of the vertices $u_{i}$, and applying \textbf{Theorem \ref{thm7}}, we arrive at
\begin{equation}
	\sum_{i=1}^{p}f^{+}(u_{i})\equiv \sum_{i=1}^{p}\Big[\Big( \sum_{v\in N(u_{i})}f(u_{i}v)\Big)(mod\hspace{1mm}p)\Big]. \label{1}
\end{equation}
From the definition of edge-graceful labeling, the vertices of the graph $G$ must be distinctly labeled from $0$ to $(p-1)$. Thus, 
\begin{equation}
	\sum_{i=1}^{p}f^{+}(u_{i}) = (0+1+2+...+(p-1)). \label{2}
\end{equation}
Substituting \eqref{2} to \eqref{1}, we will have,
\begin{equation}
	(0+1+2+...+(p-1))\equiv \sum_{i=1}^{p}\Big[\Big( \sum_{v\in N(u_{i})}f(u_{i}v)\Big)(mod\hspace{1mm}p)\Big]. \label{3}
\end{equation}
Next, observe that $\displaystyle\sum_{v\in N(u_{1})}f(u_{1}v)(mod\hspace{1mm}p)$ is the sum of the labels of the edges with endpoint $u_{1}$, modulo $p$. Similarly, $\displaystyle\sum_{v\in N(u_{2})}f(u_{2}v)(mod\hspace{1mm}p)$ is the sum of the labels of the edges with endpoint $u_{2}$, modulo $p$. This will be the case up to $\displaystyle\sum_{v\in N(u_{p})}f(u_{p}v)(mod\hspace{1mm}p)$ which is the sum of the labels of the edges with endpoint $u_{p}$, modulo $p$. Also, an edge connects two vertices in a graph. Thus,  $\displaystyle\sum_{i=1}^{p}\Big[\Big( \sum_{v\in N(u_{i})}f(u_{i}v)\Big)(mod\hspace{1mm}p)\Big]$ represents the twice the sum of the labels of all the edges, modulo $p$. Therefore,
\begin{equation}
	\sum_{i=1}^{p}\Big[\Big( \sum_{v\in N(u_{i})}f(u_{i}v)\Big)(mod\hspace{1mm}p)\Big] = 2(1+2...+q)(mod\hspace{1mm}p). \label{4}
\end{equation}
Consequently, by substituting \eqref{4} to \eqref{3} we will have,
\begin{equation*}
	(0+1+2+...+(p-1))\equiv 2(1+2...+q)(mod\hspace{1mm}p).
\end{equation*}
Now, by  {\bf Theorem \ref{thm6}}, we will have,
\begin{align*}
	2(1+2+....+q)&\equiv (0+1+2+....+(p-1))(mod\hspace{1mm}p)\\
	2\sum_{x=1}^{q}x&\equiv \Big(0+\sum_{x=1}^{p-1}x\Big)(mod\hspace{1mm}p).
\end{align*}
Since $\displaystyle\sum_{x=1}^{N}x=\frac{N(N+1)}{2}$, we will have,
\begin{align*}
	2\Big(\dfrac{q(q+1)}{2}\Big)&\equiv 0+\Big(\dfrac{(p-1)p}{2}\Big)(mod\hspace{1mm}p)\\
	(q^{2}+q)&\equiv \Big(\dfrac{(p-1)p}{2}\Big)(mod\hspace{1mm}p);
\end{align*}
which implies that,\\
\vspace{2mm}
\begin{equation*}
	p\Big|\Big(q^{2}+q-\dfrac{p(p-1)}{2}\Big).
\end{equation*} \qed

\begin{lemma}[\cite{tamang}]\label{lemma}
	\
	The quadratic Diophantine equation
	\begin{equation}
		ax^{2} + bxy + cy^{2} + dx + ey + f = 0 \label{5}
	\end{equation}
	where $a,b,c,d,e,$ and $f$ are integer coefficients, and $x$ and $y$ are the unknown variables can be reduced to
	\begin{equation}
		X^{2}-DY^{2}=N \label{6}
	\end{equation}
	where $X=Dy+E$, $Y=2ax+by+d$ and  $N=E^{2}-DF$.
\end{lemma}

Note that in this transformation, $D=b^2-4ac, E=bd-2ae, \text{and}\ F=d^2-4af$. Moreover, observe that if $X$ and $Y$ is a solution to equation (\ref{6}), then there are integers $x$ and $y$, such that
\begin{equation*}
	X=Dy+E\hspace{3mm}\rightarrow\hspace{3mm}y=\frac{X-E}{D}
\end{equation*}
\begin{equation*}
	\mathrm{and}
\end{equation*}
\begin{equation*}
	Y=2ax+by+d\hspace{3mm}\rightarrow\hspace{3mm}x=\frac{Y-by-d}{2a}
\end{equation*}
where $x$ and $y$ are the integer solutions for equation \eqref{5}.

\begin{remark}\label{remark}
	The Diophantine equations involved in this study have $c=0$. Thus, equation \eqref{6} will just be written as
	\begin{equation}
		X^{2}-(bY)^{2}=N.\label{7}
	\end{equation}
\end{remark}
Observe that we can factor the left side of equation \eqref{7} as $(X+bY)(X-bY)$. Moreover, if the pair $N_{1}$ and $N_{2}$ is a factor pair of $N$, then we will have linear system of equation,
\begin{align*}
	(X+bY)&=N_{1}\\
	(X-bY )&=N_{2}
\end{align*}
Now, solving for $X$ and $Y$, we will have,
\begin{equation*}
	X=\frac{N_{1}+N_{2}}{2}\hspace{5mm}\mathrm{and}\hspace{5mm}Y=\frac{N_{1}-N_{2}}{2b}.
\end{equation*}

\section{Results}

\label{sec3}

We are now ready to present and discuss the main result of this paper.

\begin{theorem}\label{thm12}
	Among all usual fan graphs $F_{1,n}$, only $F_{1,2}$, $F_{1,3}$, and $F_{1,11}$ are edge-graceful.
\end{theorem}

\noindent \textit{Proof.}   First, recall that a usual fan graph $F_{1,n}$ has $n+1$ vertices and $2n-1$ edges. Using {\bf{Theorem \ref{thm9}}}, which is also known as the {\bf{Lo's Theorem}}, if a graph with $p$ vertices and $q$ edges is edge graceful, then 
\begin{equation*}
	p\Big|\Big(q^{2}+q-\dfrac{p(p-1)}{2}\Big).
\end{equation*}
Thus, we substitute $p=n+1$ and $q=2n-1$ in {\bf{Theorem \ref{thm9}}}. That is,
\begin{equation*}
	(n+1)\Big|\Big((2n-1)^{2}+(2n-1)-\dfrac{(n+1)(n)}{2}\Big),
\end{equation*}
which further implies that, 
\begin{equation*}
	\Big(\dfrac{7n^{2}-5n}{2n+2}\Big)\in \mathbb{Z}.
\end{equation*}
This means,
\begin{align*}
	7n^{2}-5n &= (2n+2)k\\
	&= 2nk+2k,\hspace{5mm}\mathrm{for} \hspace{1mm} \mathrm{some} \hspace{1mm}k\in \mathbb{Z}.
\end{align*}
Thus, we solve for all the possible integer values of $n$ and $k$ in the Diophantine equation,
\begin{equation}
	7n^{2}-5n-2nk-2k  = 0 \label{8}
\end{equation}
where  $n,k\in \mathbb{Z}$.
This Diophantine equation is of the form,
\begin{equation*}
	an^{2}+bnk+ck^{2}+dn+ek+f=0,
\end{equation*}
where $a=7$,\hspace{1mm}$b=-2$,\hspace{1mm}$c=0$,\hspace{1mm}$d=-5$,\hspace{1mm}$e=-2$, and $f=0$.\\
Using the transformation stated in {\bf{Lemma \ref{lemma}}},
equation \eqref{8} can be reduced to
\begin{equation}
	{\mathbf{X^{2}-4Y^{2}=1344}} \label{9}
\end{equation}
where $Y=2an+bk+d=14n-2k-5$ and $X=Dk+E=4k+38$.\\
By factoring the left side of the equation \eqref{9}, we will have,
\begin{equation*}
	(X+2Y)(X-2Y)=1344.
\end{equation*}
Therefore, by {\bf{Lemma \ref{lemma}}} and {\bf{Remark \ref{remark}}}, we have the following:
\begin{equation*}
	X=\frac{N_{1}+N_{2}}{2},\hspace{7mm}Y=\frac{N_{1}-N_{2}}{-4},\hspace{7mm}k=\frac{X-38}{4},\hspace{7mm}n=\frac{Y+2k+5}{14},
\end{equation*}
where $(N_{1},N_{2})$ are the factor pairs of $N=1344$.
Using these information, we can now solve for the integer values of $n$ and $k$, by considering all the factor pairs of $1344$, as shown in {\bf{Table \ref{tab1}}} below.

\begin{longtable}{cccccc}
	\hline
	$N_1$&$N_2$&$X$&$Y$&$n$&$k$\\
	\hline
	1	&	1344	&	672.5	&	335.75	&	47	&	158.625	\\
2	&	672	&	337	&	167.5	&	23	&	74.75	\\
3	&	448	&	225.5	&	111.25	&	15	&	46.875	\\
4	&	336	&	170	&	83	&	{\bf\sethlcolor{yellow}\hl{11}}	&	{\bf\sethlcolor{yellow}\hl{33}}	\\
6	&	224	&	115	&	54.5	&	7	&	19.25	\\
7	&	192	&	99.5	&	46.25	&	5.857142857	&	15.375	\\
8	&	168	&	88	&	40	&	5	&	12.5	\\
12	&	112	&	62	&	25	&	{\bf\sethlcolor{yellow}\hl{3}}	&	{\bf\sethlcolor{yellow}\hl{6}}	\\
14	&	96	&	55	&	20.5	&	2.428571429	&	4.25	\\
16	&	84	&	50	&	17	&	{\bf\sethlcolor{yellow}\hl{2}}	&	{\bf\sethlcolor{yellow}\hl{3}}	\\
21	&	64	&	42.5	&	10.75	&	1.285714286	&	1.125	\\
24	&	56	&	40	&	8	&	1	&	0.5	\\
28	&	48	&	38	&	5	&	0.714285714	&	0	\\
32	&	42	&	37	&	2.5	&	0.5	&	-0.25	\\
1344	&	1	&	672.5	&	-335.75	&	-0.964285714	&	158.625	\\
672	&	2	&	337	&	-167.5	&	-0.928571429	&	74.75	\\
448	&	3	&	225.5	&	-111.25	&	-0.892857143	&	46.875	\\
336	&	4	&	170	&	-83	&	-0.857142857	&	33	\\
224	&	6	&	115	&	-54.5	&	-0.785714286	&	19.25	\\
192	&	7	&	99.5	&	-46.25	&	-0.75	&	15.375	\\
168	&	8	&	88	&	-40	&	-0.714285714	&	12.5	\\
112	&	12	&	62	&	-25	&	-0.571428571	&	6	\\
96	&	14	&	55	&	-20.5	&	-0.5	&	4.25	\\
84	&	16	&	50	&	-17	&	-0.428571429	&	3	\\
64	&	21	&	42.5	&	-10.75	&	-0.25	&	1.125	\\
56	&	24	&	40	&	-8	&	-0.142857143	&	0.5	\\
48	&	28	&	38	&	-5	&	{\bf\sethlcolor{yellow}\hl{0}}	&	{\bf\sethlcolor{yellow}\hl{0}}	\\
42	&	32	&	37	&	-2.5	&	0.142857143	&	-0.25	\\
-1	&	-1344	&	-672.5	&	-335.75	&	-49	&	-177.625	\\
-2	&	-672	&	-337	&	-167.5	&	-25	&	-93.75	\\
-3	&	-448	&	-225.5	&	-111.25	&	-17	&	-65.875	\\
-4	&	-336	&	-170	&	-83	&	{\bf\sethlcolor{yellow}\hl{-13}}	&	{\bf\sethlcolor{yellow}\hl{-52}}	\\
-6	&	-224	&	-115	&	-54.5	&	-9	&	-38.25	\\
-7	&	-192	&	-99.5	&	-46.25	&	-7.857142857	&	-34.375	\\
-8	&	-168	&	-88	&	-40	&	-7	&	-31.5	\\
-12	&	-112	&	-62	&	-25	&	{\bf\sethlcolor{yellow}\hl{-5}}	&	{\bf\sethlcolor{yellow}\hl{-25}}	\\
-14	&	-96	&	-55	&	-20.5	&	-4.428571429	&	-23.25	\\
-16	&	-84	&	-50	&	-17	&	{\bf\sethlcolor{yellow}\hl{-4}}	&	{\bf\sethlcolor{yellow}\hl{-22}}	\\
-21	&	-64	&	-42.5	&	-10.75	&	-3.285714286	&	-20.125	\\
-24	&	-56	&	-40	&	-8	&	-3	&	-19.5	\\
-28	&	-48	&	-38	&	-5	&	-2.714285714	&	-19	\\
-32	&	-42	&	-37	&	-2.5	&	-2.5	&	-18.75	\\
-1344	&	-1	&	-672.5	&	335.75	&	-1.035714286	&	-177.625	\\
-672	&	-2	&	-337	&	167.5	&	-1.071428571	&	-93.75	\\
-448	&	-3	&	-225.5	&	111.25	&	-1.107142857	&	-65.875	\\
-336	&	-4	&	-170	&	83	&	-1.142857143	&	-52	\\
-224	&	-6	&	-115	&	54.5	&	-1.214285714	&	-38.25	\\
-192	&	-7	&	-99.5	&	46.25	&	-1.25	&	-34.375	\\
-168	&	-8	&	-88	&	40	&	-1.285714286	&	-31.5	\\
-112	&	-12	&	-62	&	25	&	-1.428571429	&	-25	\\
-96	&	-14	&	-55	&	20.5	&	-1.5	&	-23.25	\\
-84	&	-16	&	-50	&	17	&	-1.571428571	&	-22	\\
-64	&	-21	&	-42.5	&	10.75	&	-1.75	&	-20.125	\\
-56	&	-24	&	-40	&	8	&	-1.857142857	&	-19.5	\\
-48	&	-28	&	-38	&	5	&	{\bf\sethlcolor{yellow}\hl{-2}}	&	{\bf\sethlcolor{yellow}\hl{-19}}	\\
-42	&	-32	&	-37	&	2.5	&	-2.142857143	&	-18.75	\\\\
	\hline
	\caption{Solutions for $7n^{2}-5n-2nk-2k  = 0 $}
	\label{tab1}
\end{longtable}

{\bf{Table \ref{tab1}}} shows that there are 8 possible solutions for the Diophantine equation $7n^{2}-5n-2nk-2k  = 0$. These are:
{\bf(n,k)=\hspace{2mm}(11,33),\hspace{2mm}(3,6),\hspace{2mm}(2,3),\hspace{2mm}(0,0),\hspace{2mm}(-13,-52),\hspace{2mm}(-5,-25),\hspace{2mm}(-4,-22) and (-2,-19)}. However, since we are only concern with possible values of $n$ for usual fan graph $F_{1,n}$, then $n$ must be a positive integer only. Hence, from our 8 solutions, only {\bf (2,3)},{\bf (3,6)}, and {\bf (11,33)} where $n=2,3$ and $11$, respectively, will satisfy this condition. Therefore, we have shown that {\bf $F_{1,2}$, $F_{1,3}$, and $F_{1,11}$ are the only usual fan graphs that satisfy the Lo's Theorem}.

\begin{figure}[h]
	\centering
	\includegraphics[width=12cm, height=10cm]{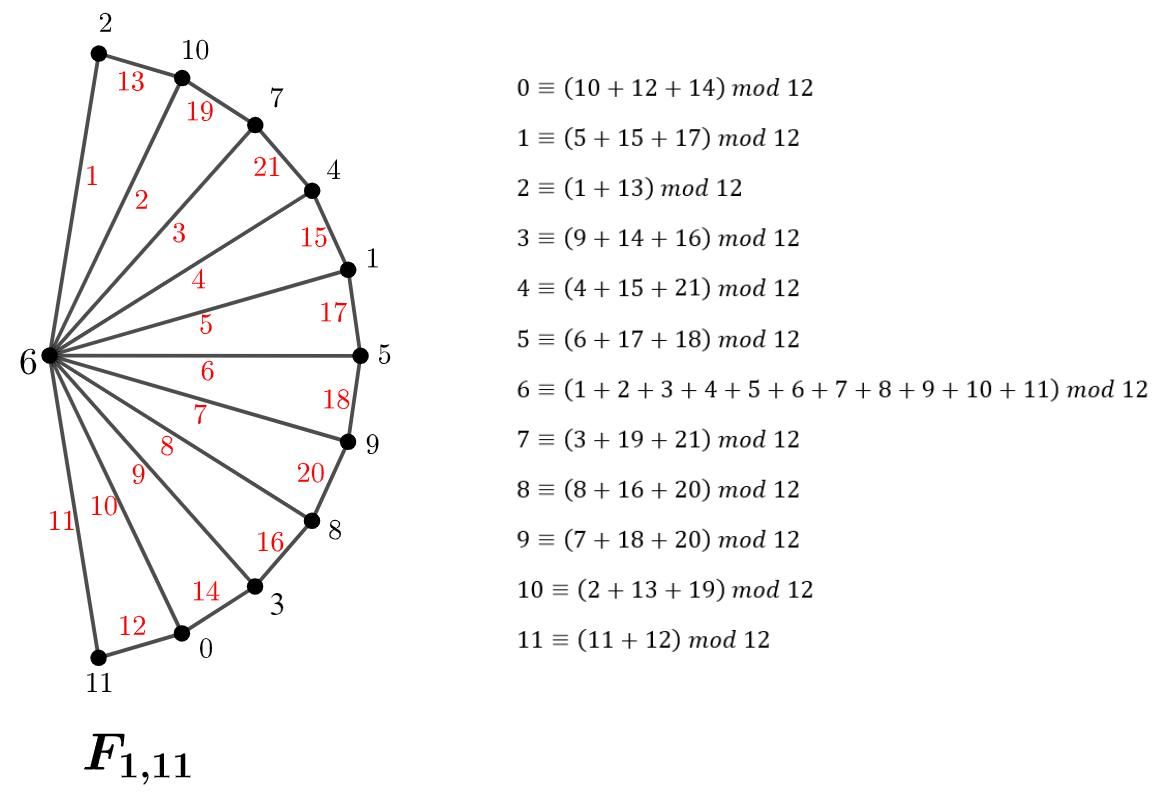}
	\caption{An Edge-graceful labeling of $F_{1,11}$}\label{fig4}
\end{figure}

\begin{figure}[h]
	\centering
	\includegraphics[width=8cm, height=4cm]{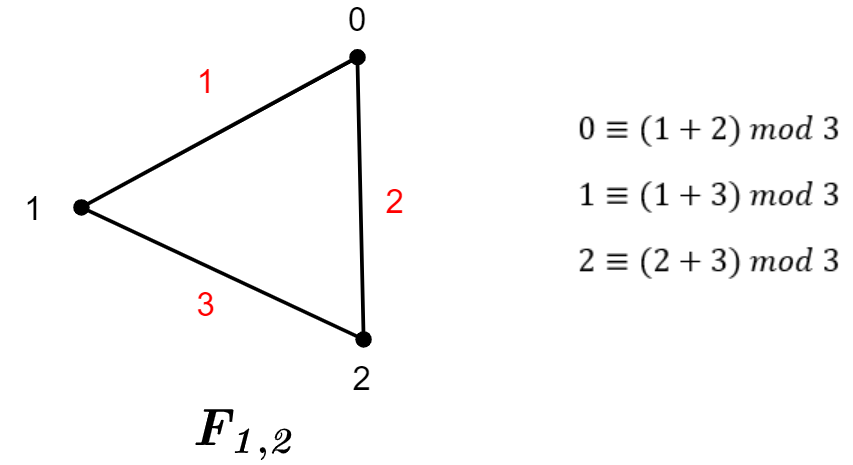}
	\caption{An Edge-graceful labeling of $F_{1,2}$}\label{fig5}
\end{figure}

\begin{figure}[h]
	\centering
	\includegraphics[width=9cm, height=5.5cm]{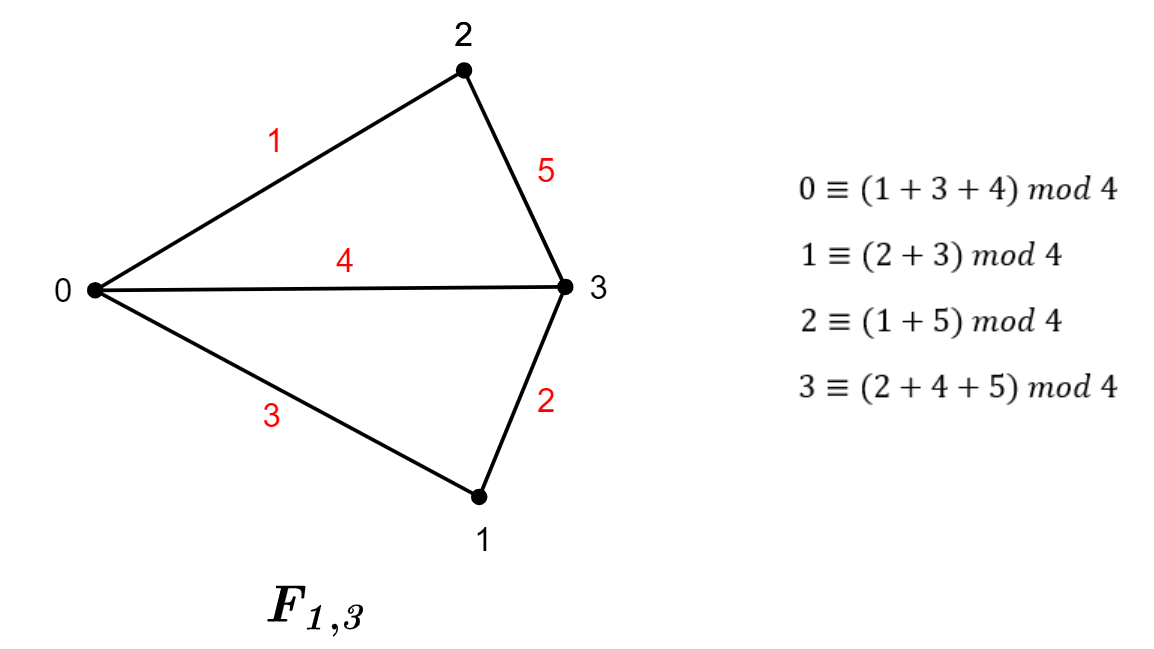}
	\caption{An Edge-graceful labeling of $F_{1,3}$}\label{fig6}
\end{figure}

To complete this proof, an actual example of an edge-graceful labeling of $F_{1,11}$, $F_{1,2}$, and $F_{1,3}$ are provided in {\bf{Figure \ref{fig4}}}, {\bf{Figure \ref{fig5}}}, and {\bf{Figure \ref{fig6}}}, respectively.

Starting with the usual fan graph $F_{1,11}$, observe that this particular graph has 12 vertices and 21 edges. Thus, by permutation, there is a total of $21!$ possible labeling for its edge. Hence, since graph $F_{1,11}$ is too large, the researchers created and used computer program on C language to generate an edge-graceful for the said graph. One of the result form the computer program is the labeling shown in {\bf{Figure \ref{fig4}}}. The notes located at the right side of the labeled-graph in {\bf{Figure \ref{fig4}}} show that it satisfies the definition of an edge-graceful labeling.

Similar computer programs were created and used to find an edge-graceful labeling for graphs $F_{1,2}$ and $F_{1,3}$. The notes located at the right side of the labeled-graph in {\bf{Figure \ref{fig5}}} and {\bf{Figure \ref{fig6}}} shows that they satisfy the definition of an edge-graceful labeling.

\newpage

Since all the needed conditions were satisfied, then we can now conclude that {\bf $F_{1,2}$, $F_{1,3}$, and $F_{1,11}$ are the only edge-graceful usual fan graphs}.\qed

\section{Conclusion and Future Works}
\label{conclusion}

In this paper, using Lo's Theorem, the concept of divisibility, solving the Diophantine equation, and a computer program, we have successfully found all edge-graceful usual fan graphs.  They are $F_{1,2}$, $F_{1,3}$, and $F_{1,11}$. The determination of edge-graceful graphs among other families of graphs, such as trees and circulant graphs, is highly recommended as a research study.  
\bigskip

{\bf Acknowledgment.} The completion of this paper would not be possible without the help and support of the following people and organizations: Central Luzon State University (CLSU), CLSU, College of Science, Department of Mathematics and Physics, Department of Science and Technology-Science Education Institute (DOST-SEI), so our heartfelt gratitude to them. We also express our gratitude to various reviewers for their comments and suggestions that helped improve the overall content of this paper.    

\bibliographystyle{amsplain}

\appendix

\section{Computer Programs Used in the Paper}
\label{app1}

This section shows all the computer programs created and used for this study. These are computer programs in C language and were used to generate actual edge-graceful labeling of needed particular graphs, in order to complete the proof of the main results  of this study.

Below, is the C program, used for generating an actual edge-graceful labeling of $F_{1,11}$\\

\begin{lstlisting}
	//usual fan graph F_{1,11} program
	
	#include <stdio.h>
	#include <stdbool.h>
	
	#define MAX_SIZE 21
	
	int Edge[MAX_SIZE];
	int Vertex[MAX_SIZE];
	int count = 1;
	int n;
	int numEdges;
	
	
	void vertices(int n){
		
		Vertex[0] = (Edge[0] + Edge[1] + Edge[2] + Edge[3] + Edge[4] + Edge[5] + Edge[6] + Edge[7] + Edge[8] + Edge[9] + Edge[10] +11) % 12;
		Vertex[1] = (Edge[10] + Edge[11] +2) % 12;
		Vertex[2] = (Edge[11] + Edge[9] + Edge[12] +3) % 12;
		Vertex[3] = (Edge[12] + Edge[8] + Edge[13] +3) % 12;
		Vertex[4] = (Edge[13] + Edge[7] + Edge[14] +3) % 12;
		Vertex[5] = (Edge[14] + Edge[6] + Edge[15] +3) % 12;
		Vertex[6] = (Edge[15] + Edge[5] + Edge[16] +3) % 12;
		Vertex[7] = (Edge[16] + Edge[4] + Edge[17] +3) % 12;
		Vertex[8] = (Edge[17] + Edge[3] + Edge[18] +3) % 12;
		Vertex[9] = (Edge[18] + Edge[2] + Edge[19] +3) % 12;
		Vertex[10] =(Edge[19] + Edge[1] + Edge[20] +3) % 12;
		Vertex[11] =(Edge[21] + +Edge[0] +2) % 12;
		
	}
	
	bool isSafe(int row, int col) {
		int i;
		for (i = 0; i < row; i++) {
			if (Edge[i] == col)
			return false;
		}
		
		return true;
		
	}
	
	
	
	void printEdge() {
		int i;
		
		for (i = 0; i < numEdges; i++) {
			printf("%d ",Edge[i] + 1);
		}
		printf("--->\n");
		for(i = 0; i <= 11; i++){
			printf("%d ", Vertex[i]);
			
		}
		printf("\n\n");
		//     }
}

void combinationofedge(int row) {
	int col;
	if (row >= numEdges) {
		vertices(n);
		
		if(Vertex[0] != Vertex[1] && Vertex[0] != Vertex[2]
		&& Vertex[0] != Vertex[3] &&Vertex[0] != Vertex[4]
		&& Vertex[0] != Vertex[5] && Vertex[0] != Vertex[6]
		&& Vertex[0] != Vertex[7] && Vertex[0] != Vertex[8]
		&& Vertex[0] != Vertex[9] && Vertex[0] != Vertex[10]
		&& Vertex[0] != Vertex[11]
		
		&& Vertex[1] != Vertex[2]
		&& Vertex[1] != Vertex[3] &&Vertex[1] != Vertex[4]
		&& Vertex[1] != Vertex[5] && Vertex[1] != Vertex[6]
		&& Vertex[1] != Vertex[7] && Vertex[1] != Vertex[8]
		&& Vertex[1] != Vertex[9] && Vertex[1] != Vertex[10]
		&& Vertex[1] != Vertex[11]
		
		&& Vertex[2] != Vertex[3] &&Vertex[2] != Vertex[4]
		&& Vertex[2] != Vertex[5] && Vertex[2] != Vertex[6]
		&& Vertex[2] != Vertex[7] && Vertex[2] != Vertex[8]
		&& Vertex[2] != Vertex[9] && Vertex[2] != Vertex[10]
		&& Vertex[2] != Vertex[11]
		
		&&Vertex[3] != Vertex[4]
		&& Vertex[3] != Vertex[5] && Vertex[3] != Vertex[6]
		&& Vertex[3] != Vertex[7] && Vertex[3] != Vertex[8]
		&& Vertex[3] != Vertex[9] && Vertex[3] != Vertex[10]
		&& Vertex[3] != Vertex[11]
		
		&& Vertex[4] != Vertex[5] && Vertex[4] != Vertex[6]
		&& Vertex[4] != Vertex[7] && Vertex[4] != Vertex[8]
		&& Vertex[4] != Vertex[9] && Vertex[4] != Vertex[10]
		&& Vertex[4] != Vertex[11]
		
		&& Vertex[5] != Vertex[6]
		&& Vertex[5] != Vertex[7] && Vertex[5] != Vertex[8]
		&& Vertex[5] != Vertex[9] && Vertex[5] != Vertex[10]
		&& Vertex[5] != Vertex[11]
		
		&& Vertex[6] != Vertex[7] && Vertex[6] != Vertex[8]
		&& Vertex[6] != Vertex[9] && Vertex[6] != Vertex[10]
		&& Vertex[6] != Vertex[11]
		
		&& Vertex[7] != Vertex[8]
		&& Vertex[7] != Vertex[9] && Vertex[7] != Vertex[10]
		&& Vertex[7] != Vertex[11]
		
		&& Vertex[8] != Vertex[9] && Vertex[8] != Vertex[10]
		&& Vertex[8] != Vertex[11]
		
		&& Vertex[9] != Vertex[10]
		&& Vertex[9] != Vertex[11]
		
		&& Vertex[10] != Vertex[11]
		
		
		){
			printf("Solution %d:\n", count++);
			printEdge();
		}
		
		return;
	}
	for (col = 0; col < numEdges; col++) {
		if (isSafe(row, col)) {
			Edge[row] = col;
			combinationofedge(row + 1);
		}
	}
}

void numberofcombinations(int n) {
	numEdges = n;
	combinationofedge(0);
}

int main() {
	
	n = 21;
	
	if (n > MAX_SIZE) {
		printf("Error: Maximum size exceeded.\n");
		return 0;
	}
	
	numberofcombinations(n);
	
	if (count == 0) {
		printf("No solutions found.\n");
	}
	
	return 0;
}
\end{lstlisting}

Below, is the C program, used for generating an actual edge-graceful labeling of $F_{1,2}$\\
\begin{lstlisting}
//usual fan graph F_{1,2} program

#include <stdio.h>
#include <stdbool.h>

#define MAX_SIZE 3

int Edge[MAX_SIZE];
int Vertex[MAX_SIZE];
int count = 1;
int n;
int numEdges;


void vertices(int n){
	
	Vertex[0] = (Edge[0] + Edge[1] + 2) % 3;
	Vertex[1] = (Edge[1] + Edge[2] + 2) % 3;
	Vertex[2] = (Edge[2] + Edge[0] + 2) % 3;
	
}

bool isSafe(int row, int col) {
	int i;
	for (i = 0; i < row; i++) {
		if (Edge[i] == col)
		return false;
	}
	
	return true;
	
}



void printEdge() {
	int i;
	
	for (i = 0; i < numEdges; i++) {
		printf("%d ",Edge[i]+1);
	}
	printf("--->\n");
	for(i = 0; i <= 2; i++){
		printf("%d ", Vertex[i]);
		
	}
	printf("\n\n");
}

void combinationofedge(int row) {
	int col;
	if (row >= numEdges) {
		vertices(n);
		
		if(Vertex[0] != Vertex[1] && Vertex[0] != Vertex[2]
		&& Vertex[1] != Vertex[2]
		
		
		){
			printf("Solution %d:\n", count++);
			printEdge();
		}
		
		return;
	}
	for (col = 0; col < numEdges; col++) {
		if (isSafe(row, col)) {
			Edge[row] = col;
			combinationofedge(row + 1);
		}
	}
}

void numberofcombinations(int n) {
	numEdges = n;
	combinationofedge(0);
}

int main() {
	
	n = 3;
	
	if (n > MAX_SIZE) {
		printf("Error: Maximum size exceeded.\n");
		return 0;
	}
	
	numberofcombinations(n);
	
	if (count == 0) {
		printf("No solutions found.\n");
	}
	
	return 0;
}
\end{lstlisting}

Below, is the C program, used for generating an actual edge-graceful labeling of $F_{1,3}$

\begin{lstlisting}
// usual fan graph F_{1,3} program

#include <stdio.h>
#include <stdbool.h>

#define MAX_SIZE 5

int Edge[MAX_SIZE];
int Vertex[MAX_SIZE];
int count = 1;
int n;
int numEdges;


void vertices(int n){
	
	Vertex[0] = (Edge[0] + Edge[1] + Edge[2] + 3) % 4;
	Vertex[1] = (Edge[2] + Edge[3] + 2) % 4;
	Vertex[2] = (Edge[1] + Edge[3] + Edge[4] + 3) % 4;
	Vertex[3] = (Edge[0] + Edge[4] + 2) % 4;
	
	
}

bool isSafe(int row, int col) {
	int i;
	for (i = 0; i < row; i++) {
		if (Edge[i] == col)
		return false;
	}
	
	return true;
	
}



void printEdge() {
	int i;
	
	for (i = 0; i < numEdges; i++) {
		printf("%d ",Edge[i] + 1);
	}
	printf("--->\n");
	for(i = 0; i <= 3; i++){
		printf("%d ", Vertex[i]);
		
	}
	printf("\n\n");
}

void combinationofedge(int row) {
	int col;
	if (row >= numEdges) {
		vertices(n);
		
		if(Vertex[0] != Vertex[1] && Vertex[0] != Vertex[2]
		&& Vertex[0] != Vertex[3] && Vertex[0] != Vertex[4]
		
		&& Vertex[1] != Vertex[2] && Vertex[1] != Vertex[3]
		&& Vertex[1] != Vertex[4]
		
		&& Vertex[2] != Vertex[3] && Vertex[2] != Vertex[4]
		
		
		){
			printf("Solution %d:\n", count++);
			printEdge();
		}
		
		return;
	}
	for (col = 0; col < numEdges; col++) {
		if (isSafe(row, col)) {
			Edge[row] = col;
			combinationofedge(row + 1);
		}
	}
}

void numberofcombinations(int n) {
	numEdges = n;
	combinationofedge(0);
}

int main() {
	
	n = 5;
	
	if (n > MAX_SIZE) {
		printf("Error: Maximum size exceeded.\n");
		return 0;
	}
	
	numberofcombinations(n);
	
	if (count == 0) {
		printf("No solutions found.\n");
	}
	
	return 0;
}
\end{lstlisting}

\end{document}